# Polymorphism of Representations of Universal Algebra

Aleks Kleyn

ABSTRACT. In this paper I consider the polymorpism of representations of universal algebra and tensor product of representations of universal algebra.

## CONTENTS

## 1. CONVENTIONS

(1) In [4], an arbitrary operation of algebra is denoted by letter $\omega$, and $\Omega$ is the set of operations of some universal algebra. Correspondingly, the universal algebra with the set of operations $\Omega$ is denoted as $\Omega$-algebra. Similar notations we see in [2] with small difference that an operation in the algebra is denoted by letter $f$ and $\mathcal{F}$ is the set of operations. I preferred first case of notations because in this case it is easier to see where I use operation.

(2) Let $A$ be $\Omega_1$-algebra. Let $B$ be $\Omega_2$-algebra. Notation

$$A \xrightarrow{*} B$$

means that there is representation of $\Omega_1$-algebra $A$ in $\Omega_2$-algebra $B$.

(3) Without a doubt, the reader may have questions, comments, objections. I will appreciate any response.

## 2. POLYMORPHISM OF REPRESENTATIONS

**Definition 2.1.** Let $A_1$, ..., $A_n$, $A$, $B_1$, ..., $B_n$, $B$ be universal algebras. Let, for any $k$, $k = 1$, ..., $n$,

$$f_k : A_k \xrightarrow{*} B_k$$

Aleks_Kleyn@MailAPS.org.
http://sites.google.com/site/AleksKleyn/.
http://arxiv.org/a/kleyn_a_1.
http://AleksKleyn.blogspot.com/.

be representation of $\Omega_1$-algebra $A_k$ in $\Omega_2$-algebra $B_k$. Let

$$f : A \dashrightarrow B$$

be representation of $\Omega_1$-algebra $A$ in $\Omega_2$-algebra $B$. The mapping

$$r : A_1 \times ... \times A_n \to A \quad R : B_1 \times ... \times B_n \to B$$

is called **polymorphism of representations** $f_1$, ..., $f_n$ into representation $f$, if, for any $k$, $k = 1$, ..., $n$,, provided that all variables except the variable $x_k \in A_k$ have given value, the mapping $(r, R)$ is a morphism of representation $f_k$ into representation $f$.

If $f_1 = ... = f_n$, then we say that the mapping $(r, R)$ is polymorphism of representation $f_1$ into representation $f$.

If $f_1 = ... = f_n = f$, then we say that the mapping $(r, R)$ is polymorphism of representation $f$. □

**Theorem 2.2.** *Let the mapping $(r, R)$ be polymorphism of representations $f_1$, ..., $f_n$ into representation $f$. The mapping $(r, R)$ satisfies to the equation*

$$(2.1) \quad R(f_1(a_1)(m_1), ..., f_n(a_n)(m_n)) = f(r(a_1, ..., a_n))(R(m_1, ..., m_n))$$

*Let $\omega_1 \in \Omega_1(p)$. For any $k$, $k = 1$, ..., $n$, the mapping $r$ satisfies to the equation*

$$(2.2) \quad \begin{aligned} &r(a_1, ..., a_{k \cdot 1}...a_{k \cdot p}\omega_1, ..., a_n) \\ &= r(a_1, ..., a_{k \cdot 1}, ..., a_n)...r(a_1, ..., a_{k \cdot p}, ..., a_n)\omega_1 \end{aligned}$$

*Let $\omega_2 \in \Omega_2(p)$. For any $k$, $k = 1$, ..., $n$, the mapping $R$ satisfies to the equation*

$$(2.3) \quad \begin{aligned} &R(m_1, ..., m_{k \cdot 1}...m_{k \cdot p}\omega_2, ..., m_n) \\ &= R(m_1, ..., m_{k \cdot 1}, ..., m_n)...R(m_1, ..., m_{k \cdot p}, ..., m_n)\omega_2 \end{aligned}$$

*Proof.* The equation (2.1) follows from the definition 2.1 and the equation [3]-(2.2.4). The equation (2.2) follows from the statement that for any $k$, $k = 1$, ..., $n$, provided that all variables except the variable $x_k \in A_k$ have given value, the mapping $r$ is homomorphism of $\Omega_1$-algebra $A_k$ into $\Omega_1$-algebra $A$. The equation (2.3) follows from the statement that for any $k$, $k = 1$, ..., $n$, provided that all variables except the variable $m_k \in B_k$ have given value, the mapping $R$ is homomorphism of $\Omega_2$-algebra $B_k$ into $\Omega_2$-algebra $B$. □

**Definition 2.3.** Let $A$, $B_1$, ..., $B_n$, $B$ be universal algebras. Let, for any $k$, $k = 1$, ..., $n$,

$$f_k : A \dashrightarrow B_k$$

be representation of $\Omega_1$-algebra $A_k$ in $\Omega_2$-algebra $B_k$. Let

$$f : A \dashrightarrow B$$

be representation of $\Omega_1$-algebra $A$ in $\Omega_2$-algebra $B$. The mapping

$$R : B_1 \times ... \times B_n \to B$$

is called **reduced polymorphism of representations** $f_1$, ..., $f_n$ into representation $f$, if, for any $k$, $k = 1$, ..., $n$, provided that all variables except the variable $x_k \in A_k$ have given value, the mapping $(\mathrm{id}, R)$ is a morphism of representation $f_k$ into representation $f$.

If $f_1 = ... = f_n$, then we say that the mapping $R$ is reduced polymorphism of representation $f_1$ into representation $f$.

If $f_1 = ... = f_n = f$, then we say that the mapping $R$ is reduced polymorphism of representation $f$. □

**Theorem 2.4.** *Let the mapping $R$ be reduced polymorphism of representations $f_1$, ..., $f_n$ into representation $f$. For any $k$, $k = 1$, ..., $n$, the mapping $R$ satisfies to the equation*

$$(2.4) \quad R(m_1, ..., f_k(a)(m_k), ..., m_n) = f(a)(R(m_1, ..., m_n))$$

*Let $\omega_2 \in \Omega_2(p)$. For any $k$, $k = 1$, ..., $n$, the mapping $R$ satisfies to the equation*

$$(2.5) \quad \begin{aligned} & R(m_1, ..., m_{k \cdot 1}...m_{k \cdot p}\omega_2, ..., m_n) \\ & = R(m_1, ..., m_{k \cdot 1}, ..., m_n)...R(m_1, ..., m_{k \cdot p}, ..., m_n)\omega_2 \end{aligned}$$

*Proof.* The equation (2.4) follows from the definition 2.3 and the equation [3]-(2.2.45). The equation (2.5) follows from the statement that for any $k$, $k = 1$, ..., $n$, provided that all variables except the variable $m_k \in B_k$ have given value, the mapping $R$ is homomorphism of $\Omega_2$-algebra $B_k$ into $\Omega_2$-algebra $B$. □

We also say that the mapping $(r, R)$ is polymorphism of representations in $\Omega_2$-algebras $B_1$, ..., $B_n$ into representation in $\Omega_2$-algebra $B$. Similarly, we say that the mapping $R$ is reduced polymorphism of representations in $\Omega_2$-algebras $B_1$, ..., $B_n$ into representation in $\Omega_2$-algebra $B$.

**Example 2.5.** Polylinear mapping of vector spaces is reduced polymorphism of vector spaces. □

Comparison of definitions 2.1 and 2.3 shows that there is a difference between these two forms of polymorphism. This is particularly evident when comparing the difference between equations (2.1) and (2.4). If we want to be able to express the reduced polymorphism of representations using polymorphism of representations, then we must require two conditions:

(1) The representation $f$ of universal algebra contains the identity transformation $\delta$. Therefore, there exists $e \in A$ such that $f(e) = \delta$. Without loss of generality, we assume that the choice of $e \in A$ does not depend on whether we consider the representation $f_1$, ..., or $f_n$.
(2) For any $k$, $k = 1$, ..., $n$,

$$(2.6) \quad r(a_1, ..., a_n) = a_k \quad a_i = e \quad i \ne k$$

Then, provided that $a_i = e$, $i \ne k$, the equation (2.1) has form

$$(2.7) \quad R(m_1, ..., f_k(a_k)(m_k), ..., m_n) = f(r(e, ..., a_k, ..., e))(R(m_1, ..., m_n))$$

It is evident that the equation (2.7) coincides with the equation (2.4).

A similar problem appears in the analysis of reduced polymorphism of representations. Consider an expression

$$(2.8) \quad R(m_1, ..., f_k(a_k)(m_k), ..., f_l(a_l)(m_l), ..., m_n)$$

We can write the equation (2.4) in the following form

$$(2.9) \quad R(m_1, ..., f_k(a) \circ m_k, ..., m_n) = f(a) \circ R(m_1, ..., m_n)$$

Using equation (2.9), we can transform the expression (2.8)

$$(2.10) \quad \begin{aligned} & R(m_1, ..., f_k(a_k) \circ m_k, ..., f_l(a_l) \circ m_l, ..., m_n) \\ & = f(a_k) \circ f(a_l) \circ R(m_1, ..., m_n) \end{aligned}$$

However, in the equation (2.10), it is not evident in what order we must consider the superposition of mappings $f(a_k)$ and $f(a_l)$. Moreover, not every $\Omega_1$-algebra $A$ has such $a$ which depends on $a_k$ and $a_l$ and satisfies to equation[1]

$$f(a) = f(a_k) \circ f(a_l)$$

The necessity of the above requirements becomes more evident, if we consider the following theorem.

**Theorem 2.6.** *Let $R$ be reduced polymorphism of representations $f_1$, ..., $f_n$ into representation $f$. Then for any $k$, $l$, $k = 1$, ..., $n$, $l = 1$, ..., $n$,*

$$(2.11) \quad \begin{aligned} & R(m_1, ..., f_l(a) \circ m_k, ..., m_l, ..., m_n) \\ & = R(m_1, ..., m_k, ..., f_l(a) \circ m_l, ..., m_n) \end{aligned}$$

*Proof.* The equation (2.11) directly follows from the equation (2.9). $\square$

Therefore, hereinafter we will require the existence of identity mapping in the representation $f$ and the existence of the mapping $r$, satisfying the equation (2.6).

## 3. Congruence

**Theorem 3.1.** *Let $N$ be equivalence on the set $A$. Let us consider category $\mathcal{A}$ whose objects are mappings*[2]

$$f_1 : A \to S_1 \quad \ker f_1 \supseteq N$$
$$f_2 : A \to S_2 \quad \ker f_2 \supseteq N$$

*We define morphism $f_1 \to f_2$ to be mapping $h : S_1 \to S_2$ making following diagram commutative*

$$\begin{array}{ccc} & & S_1 \\ & \nearrow^{f_1} & \\ A & & \downarrow h \\ & \searrow_{f_2} & \\ & & S_2 \end{array}$$

*The mapping*

$$\operatorname{nat} N : A \to A/N$$

*is universally repelling in the category $\mathcal{A}$.*[3]

---

[1] If in $\Omega_1$-algebra there is commutative product with unit, and the representation $f$ satisfies to the equation

$$f(ab) = f(a) \circ f(b)$$

then we assume

$$r(a_1, ..., a_n) = a_1...a_n$$

[2] The statement of lemma is similar to the statement on p. [1]-119.

[3] See definition of universal object of category in definition on p. [1]-57.

*Proof.* Consider diagram

$$\begin{array}{ccc} & & A/N \\ {\scriptstyle j=\mathrm{nat}\,N} \nearrow & & \big\downarrow {\scriptstyle h} \\ A & & \\ {\scriptstyle f} \searrow & & \\ & & S \end{array}$$

(3.1)  $$\ker f \supseteq N$$

From the statement (3.1) and the equation

$$j(a_1) = j(a_2)$$

it follows that

$$f(a_1) = f(a_2)$$

Therefore, we can uniquely define the mapping $h$ using the equation

$$h(j(b)) = f(b)$$

□

**Theorem 3.2.** *Let*

$$f : A \dashrightarrow B$$

*be representation of $\Omega_1$-algebra $A$ in $\Omega_2$-algebra $B$. Let $N$ be such congruence[4] on $\Omega_2$-algebra $B$ that any transformation $h \in {}^*B$ is coordinated with congruence $N$. There exists representation*

$$f_1 : A \dashrightarrow B/N$$

*of $\Omega_1$-algebra $A$ in $\Omega_2$-algebra $B/N$ and the mapping*

$$\mathrm{nat}\,N : B \to B/N$$

*is morphism of representation $f$ into the representation $f_1$*

$$\begin{array}{ccc} B & \xrightarrow{j} & B/N \quad j = \mathrm{nat}\,N \\ {\scriptstyle *f} \nwarrow & & \nearrow {\scriptstyle *f_1} \\ & A & \end{array}$$

*Proof.* We can represent any element of the set $B/N$ as $j(a)$, $a \in B$.

According to the theorem [4]-II.3.5, there exists a unique $\Omega_2$-algebra structure on the set $B/N$. If $\omega \in \Omega_2(p)$, then we define operation $\omega$ on the set $B/N$ according to the equation (3) on page [4]-59

(3.2) $$j(b_1)...j(b_p)\omega = j(b_1...b_p\omega)$$

As well as in the proof of the theorem [3]-2.2.15, we can define the representation

$$f_1 : A \dashrightarrow B/N$$

using equation

(3.3) $$f_1(a)(j(b)) = j(f(a)(b))$$

---

[4]See the definition of congruence on p. [4]-57.

We can represent the equation (3.3) using diagram

(3.4)
$$\begin{array}{ccc} B & \xrightarrow{j} & B/N \\ {\scriptstyle f(a)}\uparrow & & \uparrow{\scriptstyle f_1(a)} \\ B & \xrightarrow{j} & B/N \end{array}$$

Let $\omega \in \Omega_2(p)$. Since the mappings $f(a)$ and $j$ are homomorphisms of $\Omega_2$-algebra, then

(3.5)
$$\begin{aligned} f_1(a)(j(b_1)...j(b_p)\omega) &= f_1(a)(j(b_1...b_p\omega)) \\ &= j(f(a)(b_1...b_p\omega)) \\ &= j((f(a)(b_1))...(f(a)(b_p))\omega) \\ &= j(f(a)(b_1))...j(f(a)(b_p))\omega) \\ &= (f_1(a)(j(b_1)))...(f_1(a)(j(b_p)))\omega) \end{aligned}$$

From the equation (3.5), it follows that the mapping $f_1(a)$ is homomorphism of $\Omega_2$-algebra. From the equation (3.3) it follows that the mapping $j$ is morphism of the representation $f$ into the representation $f_1$. □

**Theorem 3.3.** *Let*
$$f : A \dashrightarrow B$$
*be representation of $\Omega_1$-algebra $A$ in $\Omega_2$-algebra $B$. Let $N$ be such congruence on $\Omega_2$-algebra $B$ that any transformation $h \in {}^*B$ is coordinated with congruence $N$. Let us consider category $\mathcal{A}$ whose objects are morphisms of representations*[5]

$$\begin{aligned} R_1 : B \to S_1 \quad &\ker R_1 \supseteq N \\ R_2 : B \to S_2 \quad &\ker R_2 \supseteq N \end{aligned}$$

*where $S_1$, $S_2$ are $\Omega_2$-algebras and*
$$g_1 : A \dashrightarrow S_1 \qquad g_2 : A \dashrightarrow S_2$$
*are representations of $\Omega_1$-algebra $A$. We define morphism $R_1 \to R_2$ to be morphism of representations $h : S_1 \to S_2$ making following diagram commutative*

$$\begin{array}{c} \xymatrix{ & & S_1 \\ A \ar@{-*>}[r]^{f} \ar@{-*>}[urr]^{g_1} \ar@{-*>}[drr]_{g_2} & B \ar[ur]^{R_1} \ar[dr]_{R_2} & \\ & & S_2 } \end{array}$$

*The morphism $\operatorname{nat} N$ of representation $f$ into representation $f_1$ (the theorem 3.2) is universally repelling in the category $\mathcal{A}$.*[6]

---

[5] The statement of lemma is similar to the statement on p. [1]-119.

[6] See definition of universal object of category in definition on p. [1]-57.

*Proof.* From the theorem 3.1, it follows that there exists and unique the mapping $h$ for which the following diagram is commutative

$$\begin{array}{c} A \xrightarrow{f}_{*} B \xrightarrow{f_1}_{*} B/N \quad j = \operatorname{nat} N \quad \ker R \supseteq N \\ \text{with } j, h, R, g \text{ maps to } S \end{array}$$

Therefore, we can uniquely define the mapping $h$ using equation

(3.6) $$h(j(b)) = R(b)$$

Let $\omega \in \Omega_2(p)$. Since mappings $R$ and $j$ are homomorphisms of $\Omega_2$-algebra, then

(3.7) $$\begin{aligned} h(j(b_1)...j(b_p)\omega) &= h(j(b_1...b_p\omega)) \\ &= R(b_1...b_p\omega) \\ &= R(b_1)...R(b_p)\omega \\ &= h(j(b_1))...h(j(b_p))\omega \end{aligned}$$

From the equation (3.7) it follows that the mapping $h$ is homomorphism of $\Omega_2$-algebra.

Since the mapping $R$ is morphism of the representation $f$ into the representation $g$, then the following equation is satisfied

(3.8) $$g(a)(R(b)) = R(f(a)(b))$$

From the equation (3.6) it follows that

(3.9) $$g(a)(h(j(b))) = g(a)(R(b))$$

From the equations (3.8), (3.9) it follows that

(3.10) $$g(a)(h(j(b))) = R(f(a)(b))$$

From the equations (3.6), (3.10) it follows that

(3.11) $$g(a)(h(j(b))) = h(j(f(a)(b)))$$

From the equations (3.3), (3.11) it follows that

(3.12) $$g(a)(h(j(b))) = h(f_1(a)(j(b)))$$

From the equation (3.12) it follows that the mapping $h$ is morphism of representation $f_1$ into the representation $g$. □

## 4. Tensor Product of Representations

**Definition 4.1.** Let $A$, $B_1$, ..., $B_n$ be universal algebras.[7] Let, for any $k$, $k = 1$, ..., $n$,

$$f_k : A \xrightarrow{\;\;*\;\;} B_k$$

---

[7] I give definition of tensor product of representations of universal algebra following to definition in [1], p. 601 - 603.

be representation of $\Omega_1$-algebra $A_k$ in $\Omega_2$-algebra $B_k$. Let us consider category $\mathcal{A}$ whose objects are reduced polymorphisms of representations $f_1$, ..., $f_n$

$$R_1 : B_1 \times ... \times B_n \longrightarrow S_1 \qquad R_2 : B_1 \times ... \times B_n \longrightarrow S_2$$

where $S_1$, $S_2$ are $\Omega_2$-algebras and

$$g_1 : A \relbar\joinrel\twoheadrightarrow S_1 \qquad g_2 : A \relbar\joinrel\twoheadrightarrow S_2$$

are representations of $\Omega_1$-algebra $A$. We define morphism $g_1 \to g_2$ to be morphisn of representations $h : S_1 \to S_2$ making following diagram commutative

$$\begin{array}{c}
\phantom{B_1 \times ... \times B_n} \xrightarrow{g_1} S_1 \\
B_1 \times ... \times B_n \phantom{\xrightarrow{g_1}} \phantom{S_1} \downarrow h \\
\phantom{B_1 \times ... \times B_n} \xrightarrow[g_2]{} S_2
\end{array}$$

Universal object $B_1 \otimes ... \otimes B_n$ of category $\mathcal{A}$ is called **tensor product of representations** $B_1$, ..., $B_n$. □

**Definition 4.2.** Tensor product

$$B^{\otimes n} = B_1 \otimes ... \otimes B_n \quad B_1 = ... = B_n = B$$

is called **tensor power of representation** $B$. □

**Theorem 4.3.** *There exists tensor product of representations.*

*Proof.* Let

$$f : A \relbar\joinrel\twoheadrightarrow M$$

be representation of $\Omega_1$-algebra $A$ generated by product $B_1 \times ... \times B_n$ of $\Omega_2$-algebras $B_1$, ..., $B_n$ (the teorem [3]-3.1.4). Injection

$$i : B_1 \times ... \times B_n \longrightarrow M$$

is defined according to rule

(4.1) $$i \circ (b_1, ..., b_n) = (b_1, ..., b_n)$$

Let $N$ be equivalence generated by following equations

(4.2) $$(b_1, ..., b_{i \cdot 1} ... b_{i \cdot p} \omega, ..., b_n) = (b_1, ..., b_{i \cdot 1}, ..., b_n)...(b_1, ..., b_{i \cdot p}, ..., b_n)\omega$$

(4.3) $$(b_1, ..., f_i(a)(b_i), ..., b_n) = f(a)((b_1, ..., b_i, ..., b_n))$$

$$b_k \in B_k \quad k = 1, ..., n \quad b_{i \cdot 1}, ..., b_{i \cdot p} \in B_i \quad \omega \in \Omega_2(p) \quad a \in A$$

- Let us prove the following lemma

    *Lemma 4.4.* For any $c \in A$, endomorphism $f(c)$ of $\Omega_2$-algebra $M$ is coordinated with equivalence $N$.

- Let $\omega \in \Omega_2(p)$. From the equation (4.3), it follows that

(4.4) $\quad f(c)((b_1, ..., b_{i\cdot 1}...b_{i\cdot p}\omega, ..., b_n)) = (b_1, ..., f_i(c)(b_{i\cdot 1}...b_{i\cdot p}\omega), ..., b_n)$

Since $f_i(c)$ is endomorphism of $\Omega_2$-algebra $B_i$, then from the equation (4.4), it follows that

(4.5) $\quad f(c)((b_1, ..., b_{i\cdot 1}...b_{i\cdot p}\omega, ..., b_n)) = (b_1, ..., f_i(c)(b_{i\cdot 1})...f_i(c)(b_{i\cdot p})\omega, ..., b_n)$

From equations (4.5), (4.2), it follows that

(4.6) $\quad \begin{aligned} & f(c)((b_1, ..., b_{i\cdot 1}...b_{i\cdot p}\omega, ..., b_n)) \\ & = (b_1, ..., f_i(c)(b_{i\cdot 1}), ..., b_n)...(b_1, ..., f_i(c)(b_{i\cdot p}), ..., b_n)\omega \end{aligned}$

From equations (4.6), (4.3), it follows that

(4.7) $\quad \begin{aligned} & f(c)((b_1, ..., b_{i\cdot 1}...b_{i\cdot p}\omega, ..., b_n)) \\ & = f(c)((b_1, ..., b_{i\cdot 1}, ..., b_n))...f(c)((b_1, ..., b_{i\cdot p}, ..., b_n))\omega \end{aligned}$

Since $f_i(c)$ is endomorphism of $\Omega_2$-algebra $B_i$, then from the equation (4.7), it follows that

(4.8) $\quad \begin{aligned} & f(c)((b_1, ..., b_{i\cdot 1}...b_{i\cdot p}\omega, ..., b_n)) \\ & = f(c)((b_1, ..., b_{i\cdot 1}, ..., b_n))...(b_1, ..., b_{i\cdot p}, ..., b_n)\omega) \end{aligned}$

Therefore, we proved the following statement.

*Lemma* 4.5. The equation (4.8) follows from the equation (4.2).

- The following statement follows from the equation (4.3).

*Lemma* 4.6. The equation

(4.9) $\quad f(c)((b_1, ..., f_i(a)(b_i), ..., b_n)) = f(c)(f(a)((b_1, ..., b_i, ..., b_n)))$

follows from the equation (4.3).

- The lemma 4.4 follows from lemmas 4.5, 4.6 and definition [3]-2.2.13.

From the lemma 4.4 and the theorem [3]-2.2.14, it follows that $\Omega_1$-algebra is defined on the set $^*M/N$. Consider diagram

$$\begin{array}{ccc} M/N & \xrightarrow{F(a)} & M/N \quad j = \text{nat } N \\ \uparrow \nearrow^{F} & & \uparrow \\ A \quad j & & j \\ \downarrow \searrow_{f} & & \uparrow \\ M & \xrightarrow{f(a)} & M \end{array}$$

According to lemma 4.4, from the condition

$$j(b_1) = j(b_2)$$

it follows that
$$j(f(a)(b_1)) = j(f(a)(b_2))$$
Therefore, transformation $F(a)$ is well defined and
$$(4.10) \quad F(a) \circ j = j \circ f(a)$$
If $\omega \in \Omega_1(p)$, then we assume
$$(F(a_1)...F(a_p)\omega)(J(b)) = J((f(a_1)...f(a_p)\omega)(b))$$
Therefore, mapping $F$ is representations of $\Omega_1$-algebra $A$. From (4.10) it follows that $(\text{id}, j)$ is morphism of representations $f$ and $F$.

Consider commutative diagram

(4.11)
$$\begin{array}{c} B_1 \times ... \times B_n \xrightarrow{i} M \xrightarrow{j} M/N \\ \text{with } g_1 : B_1 \times ... \times B_n \to M/N \end{array}$$

From equations (4.1), (4.2), (4.3), it follows that

(4.12)
$$g_1((b_1, ..., b_{i \cdot 1}...b_{i \cdot p}\omega, ..., b_n))$$
$$= g_1((b_1, ..., b_{i \cdot 1}, ..., b_n))...g_1((b_1, ..., b_{i \cdot p}, ..., b_n))\omega$$

(4.13)
$$g_1((b_1, ..., f_i(a)(b_i), ..., b_n))$$
$$= f(a)(g_1((b_1, ..., b_i, ..., b_n)))$$

From equations (4.12) and (4.13) it follows that map $g_1$ is reduced polymorphism of representations $f_1, ..., f_n$.

Since $B_1 \times ... \times B_n$ is the basis of representation $M$ of $\Omega_1$ algebra $A$, then, according to the theorem [3]-3.2.7, for any representation
$$A \xrightarrow{*} V$$
and any reduced polymorphism
$$g_2 : B_1 \times ... \times B_n \longrightarrow V$$
there exists a unique morphism of representations $k : M \to V$, for which following diagram is commutative

(4.14)
$$\begin{array}{c} B_1 \times ... \times B_n \xrightarrow{i} M \\ \searrow_{g_2} \quad \downarrow_k \\ V \end{array}$$

Since $g_2$ is reduced polymorphism, then $\ker k \supseteq N$.

According to the theorem 3.3, map $j$ is universal in the category of morphisms of representation $f$ whose kernel contains $N$. Therefore, we have morphism of representations
$$h : M/N \to V$$

which makes the following diagram commutative

(4.15)
$$\begin{array}{c} M/N \\ {}^{j}\nearrow \quad \downarrow h \\ M \quad\quad \\ {}_{k}\searrow \quad \\ V \end{array}$$

We join diagrams (4.11), (4.14), (4.15), and get commutative diagram

$$\begin{array}{c} \xrightarrow{g_1} M/N \\ B_1\times...\times B_n \xrightarrow{i} M \quad \downarrow h \\ \xrightarrow{g_2} \searrow^k \\ V \end{array}$$

Since $\operatorname{Im} g_1$ generates $M/N$, than map $h$ is uniquely determined. $\square$

According to proof of theorem 4.3

$$B_1 \otimes ... \otimes B_n = M/N$$

If $d_i \in A_i$, we write

(4.16) $$j \circ (d_1, ..., d_n) = d_1 \otimes ... \otimes d_n$$

**Theorem 4.7.** *Let $B_1$, ..., $B_n$ be $\Omega_2$-algebras. Let*

$$f : B_1 \times ... \times B_n \to B_1 \otimes ... \otimes B_n$$

*be reduced polymorphism defined by equation*

(4.17) $$f \circ (b_1, ..., b_n) = b_1 \otimes ... \otimes b_n$$

*Let*

$$g : B_1 \times ... \times B_n \to V$$

*be reduced polymorphism into $\Omega$-algebra $V$. There exists morphism of representations*

$$h : B_1 \otimes ... \otimes B_n \to V$$

*such that the diagram*

$$\begin{array}{c} B_1 \otimes ... \otimes B_n \\ {}^{f}\nearrow \quad \downarrow h \\ B_1 \times ... \times B_n \quad\quad \\ {}_{g}\searrow \quad \\ V \end{array}$$

*is commutative.*

*Proof.* Equation (4.17) follows from equations (4.1) and (4.16). An existence of the mapping $h$ follows from the definition 4.1 and constructions made in the proof of the theorem 4.3. □

We can write equations (4.12) and (4.13) as

$$
\begin{aligned}
&b_1 \otimes ... \otimes (b_{i \cdot 1}...b_{i \cdot p}\omega) \otimes ... \otimes b_n \\
&= (b_1 \otimes ... \otimes b_{i \cdot 1} \otimes ... \otimes b_n)...(b_1 \otimes ... \otimes b_{i \cdot p} \otimes ... \otimes b_n)\omega
\end{aligned}
\tag{4.18}
$$

$$
b_1 \otimes ... \otimes (f_i(a)(b_i)) \otimes ... \otimes b_n = f(a)(b_1 \otimes ... \otimes b_i \otimes ... \otimes b_n) \tag{4.19}
$$

$$b_k \in B_k \quad k = 1,...,n \quad b_{i \cdot 1},...,b_{i \cdot p} \in B_i \quad \omega \in \Omega_2(p) \quad a \in A$$

# 6. Index

# 7. Special Symbols and Notations

# Полиморфим представлений универсальной алгебры

Александр Клейн

Аннотация. В статье я рассматриваю полиморфизм представлений универсальной алгебры и тензорное произведение представлений универсальной алгебры.

Содержание

## 1. Соглашения

(1) В [4] произвольная операция алгебры обозначена буквой $\omega$, и $\Omega$ - множество операций некоторой универсальной алгебры. Соответственно, универсальная алгебра с множеством операций $\Omega$ обозначается $\Omega$-алгебра. Аналогичные обозначения мы видим в [2] с той небольшой разницей, что операция в алгебре обозначена буквой $f$ и $\mathcal{F}$ - множество операций. Я выбрал первый вариант обозначений, так как в этом случае легче видно, где я использую операцию.

(2) Пусть $A$ - $\Omega_1$-алгебра. Пусть $B$ - $\Omega_2$-алгебра. Запись

$$A \xrightarrow{*} B$$

означает, что определено представление $\Omega_1$-алгебры $A$ в $\Omega_2$-алгебре $B$.

(3) Без сомнения, у читателя могут быть вопросы, замечания, возражения. Я буду признателен любому отзыву.

## 2. Полиморфизм представлений

**Определение 2.1.** Пусть $A_1$, ..., $A_n$, $A$, $B_1$, ..., $B_n$, $B$ - универсальные алгебры. Пусть для любого $k$, $k = 1, ..., n$,

$$f_k : A_k \xrightarrow{*} B_k$$

Aleks_Kleyn@MailAPS.org.
http://sites.google.com/site/AleksKleyn/.
http://arxiv.org/a/kleyn_a_1.
http://AleksKleyn.blogspot.com/.

представление $\Omega_1$-алгебры $A_k$ в $\Omega_2$-алгебре $B_k$. Пусть

$$f : A \xrightarrow{\ *\ } B$$

представление $\Omega_1$-алгебры $A$ в $\Omega_2$-алгебре $B$. Отображение

$$r : A_1 \times ... \times A_n \to A \quad R : B_1 \times ... \times B_n \to B$$

называется **полиморфизмом представлений** $f_1, ..., f_n$ в представление $f$, если для любого $k$, $k = 1, ..., n$, при условии, что все переменные кроме переменной $x_k \in A_k$ имеют заданное значение, отображение $(r, R)$ является морфизмом представления $f_k$ в представление $f$.

Если $f_1 = ... = f_n$, то мы будем говорить, что отображение $(r, R)$ является полиморфизмом представления $f_1$ в представление $f$.

Если $f_1 = ... = f_n = f$, то мы будем говорить, что отображение $(r, R)$ является полиморфизмом представления $f$. $\square$

**Теорема 2.2.** *Пусть отображение $(r, R)$ является полиморфизмом представлений $f_1, ..., f_n$ в представление $f$. Отображение $(r, R)$ удовлетворяет равенству*

(2.1) $\quad R(f_1(a_1)(m_1), ..., f_n(a_n)(m_n)) = f(r(a_1, ..., a_n))(R(m_1, ..., m_n))$

*Пусть $\omega_1 \in \Omega_1(p)$. Для любого $k$, $k = 1, ..., n$, отображение $r$ удовлетворяет равенству*

(2.2) $\quad\begin{aligned}&r(a_1, ..., a_{k\cdot 1}...a_{k\cdot p}\omega_1, ..., a_n) \\ &= r(a_1, ..., a_{k\cdot 1}, ..., a_n)...r(a_1, ..., a_{k\cdot p}, ..., a_n)\omega_1\end{aligned}$

*Пусть $\omega_2 \in \Omega_2(p)$. Для любого $k$, $k = 1, ..., n$, отображение $R$ удовлетворяет равенству*

(2.3) $\quad\begin{aligned}&R(m_1, ..., m_{k\cdot 1}...m_{k\cdot p}\omega_2, ..., m_n) \\ &= R(m_1, ..., m_{k\cdot 1}, ..., m_n)...R(m_1, ..., m_{k\cdot p}, ..., m_n)\omega_2\end{aligned}$

*Доказательство.* Равенство (2.1) следует из определения 2.1 и равенства [3]-(2.2.4). Равенство (2.2) следует из утверждения, что для любого $k$, $k = 1, ..., n$, при условии, что все переменные кроме переменной $x_k \in A_k$ имеют заданное значение, отображение $r$ является гомоморфизмом $\Omega_1$-алгебры $A_k$ в $\Omega_1$-алгебру $A$. Равенство (2.3) следует из утверждения, что для любого $k$, $k = 1, ..., n$, при условии, что все переменные кроме переменной $m_k \in B_k$ имеют заданное значение, отображение $R$ является гомоморфизмом $\Omega_2$-алгебры $B_k$ в $\Omega_2$-алгебру $B$. $\square$

**Определение 2.3.** Пусть $A$, $B_1$, ..., $B_n$, $B$ - универсальные алгебры. Пусть для любого $k$, $k = 1, ..., n$,

$$f_k : A \xrightarrow{\ *\ } B_k$$

представление $\Omega_1$-алгебры $A_k$ в $\Omega_2$-алгебре $B_k$. Пусть

$$f : A \xrightarrow{\ *\ } B$$

представление $\Omega_1$-алгебры $A$ в $\Omega_2$-алгебре $B$. Отображение

$$R : B_1 \times ... \times B_n \to B$$

называется **приведенным полиморфизмом представлений** $f_1$, ..., $f_n$ в представление $f$, если для любого $k$, $k = 1$, ..., $n$, при условии, что все переменные кроме переменной $x_k \in A_k$ имеют заданное значение, отображение $(\mathrm{id}, R)$ является морфизмом представления $f_k$ в представление $f$.

Если $f_1 = ... = f_n$, то мы будем говорить, что отображение $R$ является приведенным полиморфизмом представления $f_1$ в представление $f$.

Если $f_1 = ... = f_n = f$, то мы будем говорить, что отображение $R$ является приведенным полиморфизмом представления $f$. $\square$

**Теорема 2.4.** *Пусть отображение $R$ является приведенным полиморфизмом представлений $f_1$, ..., $f_n$ в представление $f$. Для любого $k$, $k = 1$, ..., $n$, отображение $R$ удовлетворяет равенству*

(2.4) $$R(m_1, ..., f_k(a)(m_k), ..., m_n) = f(a)(R(m_1, ..., m_n))$$

*Пусть $\omega_2 \in \Omega_2(p)$. Для любого $k$, $k = 1$, ..., $n$, отображение $R$ удовлетворяет равенству*

(2.5) $$\begin{aligned}&R(m_1, ..., m_{k\cdot 1}...m_{k\cdot p}\omega_2, ..., m_n) \\ &= R(m_1, ..., m_{k\cdot 1}, ..., m_n)...R(m_1, ..., m_{k\cdot p}, ..., m_n)\omega_2\end{aligned}$$

*Доказательство.* Равенство (2.4) следует из определения 2.3 и равенства [3]-(2.2.45). Равенство (2.5) следует из утверждения, что для любого $k$, $k = 1$, ..., $n$, при условии, что все переменные кроме переменной $m_k \in B_k$ имеют заданное значение, отображение $R$ является гомоморфизмом $\Omega_2$-алгебры $B_k$ в $\Omega_2$-алгебру $B$. $\square$

Мы также будем говорить, что отображение $(r, R)$ является полиморфизмом представлений в $\Omega_2$-алгебрах $B_1$, ..., $B_n$ в представление в $\Omega_2$-алгебре $B$. Аналогично, мы будем говорить, что отображение $R$ является приведенным полиморфизмом представлений в $\Omega_2$-алгебрах $B_1$, ..., $B_n$ в представление в $\Omega_2$-алгебре $B$.

**Пример 2.5.** Полилинейное отображение векторных пространств является приведенным полиморфизмом векторных пространств. $\square$

Сравнение определений 2.1 и 2.3 показывает, что существует разница между этими двумя формами полиморфизма. Особенно хорошо это различие видно при сравнении равенств (2.1) и (2.4). Если мы хотим иметь возможность выразить приведенный полиморфизм представлений через полиморфизм представлений, то мы должны потребовать, два условия:

(1) Представление $f$ универсальной алгебры содержит тождественное преобразование $\delta$. Следовательно, существует $e \in A$ такой, что $f(e) = \delta$. Не нарушая общности, мы положим, что выбор $e \in A$ не зависит от того, какое представление $f_1$, ..., $f_n$ мы рассматриваем.
(2) Для любого $k$, $k = 1$, ..., $n$,

(2.6) $$r(a_1, ..., a_n) = a_k \quad a_i = e \quad i \neq k$$

Тогда, при условии $a_i = e$, $i \neq k$, равенство (2.1) имеет вид

(2.7) $$R(m_1, ..., f_k(a_k)(m_k), ..., m_n) = f(r(e, ..., a_k, ..., e))(R(m_1, ..., m_n))$$

Очевидно, что равенство (2.7) совпадает с равенством (2.4).

Похожая задача появляется при анализе приведенного полиморфизма представлений. Рассмотрим выражение

$$(2.8) \qquad R(m_1, ..., f_k(a_k)(m_k), ..., f_l(a_l)(m_l), ..., m_n)$$

Мы можем записать равенство (2.4) в виде

$$(2.9) \qquad R(m_1, ..., f_k(a) \circ m_k, ..., m_n) = f(a) \circ R(m_1, ..., m_n)$$

Пользуясь равенством (2.9), мы можем преобразовать выражение (2.8)

$$(2.10) \qquad \begin{aligned} &R(m_1, ..., f_k(a_k) \circ m_k, ..., f_l(a_l) \circ m_l, ..., m_n) \\ &= f(a_k) \circ f(a_l) \circ R(m_1, ..., m_n) \end{aligned}$$

Однко в равенстве (2.10) не очевидно, в каком порядке мы должны рассматривать суперпозицию отображений $f(a_k)$ и $f(a_l)$. Кроме того, не всякая $\Omega_1$-алгебра $A$ имеет такой $a$ (зависящий от $a_k$ и $a_l$), что[1]

$$f(a) = f(a_k) \circ f(a_l)$$

Необходимость указанных выше требований становится более очевидной, если мы рассмотрим следующую теорему.

**Теорема 2.6.** *Пусть $R$ - приведенный полиморфизм представлений $f_1$, ..., $f_n$ в представление $f$. Тогда для любых $k$, $l$, $k = 1, ..., n$, $l = 1, ..., n$,*

$$(2.11) \qquad \begin{aligned} &R(m_1, ..., f_l(a) \circ m_k, ..., m_l, ..., m_n) \\ &= R(m_1, ..., m_k, ..., f_l(a) \circ m_l, ..., m_n) \end{aligned}$$

*Доказательство.* Равенство (2.11) непосредственно следует из равенства (2.9). □

Поэтому в дальнейшем мы будем требовать существование тождественного отображения в представлении $f$ и существование отображения $r$, удовлетворяющего равенству (2.6).

## 3. Конгруэнция

**Теорема 3.1.** *Пусть $N$ - отношение эквивалентности на множестве $A$. Рассмотрим категорию $\mathcal{A}$ объектами которой являются отображения*[2]

$$\begin{aligned} f_1 &: A \to S_1 \quad \ker f_1 \supseteq N \\ f_2 &: A \to S_2 \quad \ker f_2 \supseteq N \end{aligned}$$

---

[1] Если в $\Omega_1$-алгебре определенно коммутативное произведение с единицей, а представление $f$ удовлетворяет равенству

$$f(ab) = f(a) \circ f(b)$$

то мы положим

$$r(a_1, ..., a_n) = a_1...a_n$$

[2] Утверждение леммы аналогично утверждению на с. [1]-94.

Мы определим морфизм $f_1 \to f_2$ как отображение $h : S_1 \to S_2$, для которого коммутативна диаграмма

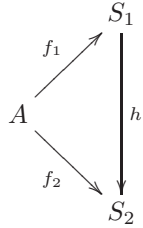

Отображение является универсально отталкивающим в категории $\mathcal{A}$.[3]

*Доказательство.* Рассмотрим диаграмму

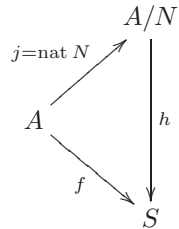

(3.1) $$\ker f \supseteq N$$

Из утверждения (3.1) и равенства

$$j(a_1) = j(a_2)$$

следует

$$f(a_1) = f(a_2)$$

Следовательно, мы можем однозначно определить отображение $h$ с помощью равенства

$$h(j(b)) = f(b)$$

$\square$

**Теорема 3.2.** *Пусть*

$$f : A \dashrightarrow B$$

*представление $\Omega_1$-алгебры $A$ в $\Omega_2$-алгебре $B$. Пусть $N$ - такая конгруэнция[4] на $\Omega_2$-алгебре $B$, что любое преобразование $h \in {}^*B$ согласованно с конгруэнцией $N$. Существует представление*

$$f_1 : A \dashrightarrow B/N$$

*$\Omega_1$-алгебры $A$ в $\Omega_2$-алгебре $B/N$ и отображение*

$$\operatorname{nat} N : B \to B/N$$

---

[3]Определение универсального объекта смотри в определении на с. [1]-47.

[4]Смотри определение конгруэнции на с. [4]-71.

*является морфизмом представления $f$ в представление $f_1$*

$$\begin{array}{c} B \xrightarrow{\ j\ } B/N \quad j = \operatorname{nat} N \\ {}_{f}{\nwarrow}^{*} \quad {}^{*}{\nearrow}_{f_1} \\ A \end{array}$$

*Доказательство.* Любой элемент множества $B/N$ мы можем представить в виде $j(a)$, $a \in B$.

Согласно теореме [4]-II.3.5, мы можем определить единственную структуру $\Omega_2$-алгебры на множестве $B/N$. Если $\omega \in \Omega_2(p)$, то мы определим операцию $\omega$ на множестве $B/N$ согласно равенству (3) на странице [4]-73

(3.2) $$j(b_1)...j(b_p)\omega = j(b_1...b_p\omega)$$

Также как в доказательстве теоремы [3]-2.2.15, мы можем определить представление
$$f_1 : A \longrightarrow^* B/N$$
с помощью равенства

(3.3) $$f_1(a)(j(b)) = j(f(a)(b))$$

Равенство (3.3) можно представить с помощью диаграммы

(3.4) $$\begin{array}{ccc} B & \xrightarrow{\ j\ } & B/N \\ {\scriptstyle f(a)}\uparrow & & \uparrow {\scriptstyle f_1(a)} \\ B & \xrightarrow{\ j\ } & B/N \end{array}$$

Пусть $\omega \in \Omega_2(p)$. Так как отображения $f(a)$ и $j$ являются гомоморфизмами $\Omega_2$-алгебры, то

(3.5) $$\begin{aligned} f_1(a)(j(b_1)...j(b_p)\omega) &= f_1(a)(j(b_1...b_p\omega)) \\ &= j(f(a)(b_1...b_p\omega)) \\ &= j((f(a)(b_1))...(f(a)(b_p))\omega) \\ &= j(f(a)(b_1))...j(f(a)(b_p))\omega \\ &= (f_1(a)(j(b_1)))...(f_1(a)(j(b_p)))\omega \end{aligned}$$

Из равенства (3.5) следует, что отображение $f_1(a)$ является гомоморфизмом $\Omega_2$-алгебры. Из равенства (3.3) следует, что отображение $j$ является морфизмом представления $f$ в представление $f_1$. $\square$

**Теорема 3.3.** *Пусть*
$$f : A \longrightarrow^* B$$
*представление $\Omega_1$-алгебры $A$ в $\Omega_2$-алгебре $B$. Пусть $N$ - такая конгруэнция на $\Omega_2$-алгебре $B$, что любое преобразование $h \in {}^*B$ согласованно с конгруэнцией $N$. Рассмотрим категорию $\mathcal{A}$ объектами которой являются морфизмы представлений*[5]

$$\begin{aligned} R_1 : B \to S_1 &\quad \ker R_1 \supseteq N \\ R_2 : B \to S_2 &\quad \ker R_2 \supseteq N \end{aligned}$$

---

[5]Утверждение леммы аналогично утверждению на с. [1]-94.

*где $S_1$, $S_2$ - $\Omega_2$-алгебры и*

$$g_1 : A \dashrightarrow S_1 \qquad g_2 : A \dashrightarrow S_2$$

*представления $\Omega_1$-алгебры $A$. Мы определим морфизм $R_1 \to R_2$ как морфизм представлений $h : S_1 \to S_2$, для которого коммутативна диаграмма*

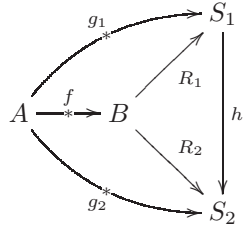

*Морфизм $\operatorname{nat} N$ представления $f$ в представление $f_1$ (теорема 3.2) является универсально отталкивающим в категории $\mathcal{A}$.*[6]

*Доказательство.* Существование и единственность отображения $h$, для которого коммутативна диаграмма

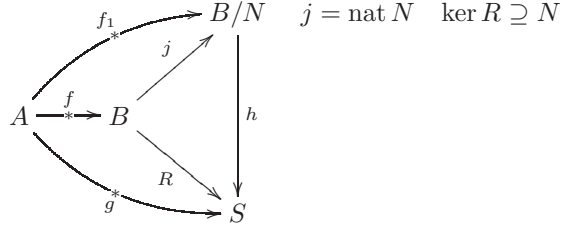

следует из теоремы 3.1. Следовательно, мы можем однозначно определить отображение $h$ с помощью равенства

$$(3.6) \qquad h(j(b)) = R(b)$$

Пусть $\omega \in \Omega_2(p)$. Так как отображения $R$ и $j$ являются гомоморфизмами $\Omega_2$-алгебры, то

$$(3.7) \quad \begin{aligned} h(j(b_1)...j(b_p)\omega) &= h(j(b_1...b_p\omega)) \\ &= R(b_1...b_p\omega) \\ &= R(b_1)...R(b_p)\omega \\ &= h(j(b_1))...h(j(b_p))\omega \end{aligned}$$

следует, что отображение $h$ является гомоморфизмом $\Omega_2$-алгебры.

Так как отображение $R$ является морфизмом представления $f$ в представление $g$, то верно равенство

$$(3.8) \qquad g(a)(R(b)) = R(f(a)(b))$$

Из равенства (3.6) следует

$$(3.9) \qquad g(a)(h(j(b))) = g(a)(R(b))$$

Из равенств (3.8), (3.9) следует

$$(3.10) \qquad g(a)(h(j(b))) = R(f(a)(b))$$

---

[6]Определение универсального объекта смотри в определении на с. [1]-47.

Из равенств (3.6), (3.10) следует

(3.11) $$g(a)(h(j(b))) = h(j(f(a)(b)))$$

Из равенств (3.3), (3.11) следует

(3.12) $$g(a)(h(j(b))) = h(f_1(a)(j(b)))$$

Из равенства (3.12) следует, что отображение $h$ является морфизмом представления $f_1$ в представление $g$. □

## 4. Тензорное произведение представлений

**Определение 4.1.** Пусть $A$, $B_1$, ..., $B_n$ - универсальные алгебры.[7] Пусть для любого $k$, $k = 1, ..., n$,

$$f_k : A \dashrightarrow B_k$$

представление $\Omega_1$-алгебры $A_k$ в $\Omega_2$-алгебре $B_k$. Рассмотрим категорию $\mathcal{A}$ объектами которой являются приведенные полиморфизмы представлений $f_1$, ..., $f_n$

$$R_1 : B_1 \times ... \times B_n \longrightarrow S_1 \qquad R_2 : B_1 \times ... \times B_n \longrightarrow S_2$$

где $S_1$, $S_2$ - $\Omega_2$-алгебры и

$$g_1 : A \dashrightarrow S_1 \quad g_2 : A \dashrightarrow S_2$$

представления $\Omega_1$-алгебры $A$. Мы определим морфизм $g_1 \to g_2$ как морфизм представлений $h : S_1 \to S_2$, для которого коммутативна диаграмма

$$\begin{array}{c} & & S_1 \\ & \nearrow_{g_1} & \\ B_1 \times ... \times B_n & & \downarrow h \\ & \searrow_{g_2} & \\ & & S_2 \end{array}$$

Универсальный объект $B_1 \otimes ... \otimes B_n$ категории $\mathcal{A}$ называется **тензорным произведением представлений** $B_1$, **...**, $B_n$. □

**Определение 4.2.** Тензорное произведение

$$B^{\otimes n} = B_1 \otimes ... \otimes B_n \quad B_1 = ... = B_n = B$$

называется **тензорной степенью представления** $B$. □

**Теорема 4.3.** *Тензорное произведение представлений существует.*

*Доказательство.* Пусть

$$f : A \dashrightarrow M$$

представление $\Omega_1$-алгебры $A$, порождённое произведением $B_1 \times ... \times B_n$ $\Omega_2$-алгебр $B_1$, ..., $B_n$ (теорема [3]-3.1.4). Инъекция

$$i : B_1 \times ... \times B_n \longrightarrow M$$

---

[7]Я определяю тензорное произведение представлений универсальной алгебры по аналогии с определением в [1], с. 456 - 458.

определена по правилу

(4.1) $$i \circ (b_1, ..., b_n) = (b_1, ..., b_n)$$

Пусть $N$ - отношение эквивалентности, порождённое равенствами

(4.2) $$(b_1, ..., b_{i \cdot 1}...b_{i \cdot p}\omega, ..., b_n) = (b_1, ..., b_{i \cdot 1}, ..., b_n)...(b_1, ..., b_{i \cdot p}, ..., b_n)\omega$$

(4.3) $$(b_1, ..., f_i(a)(b_i), ..., b_n) = f(a)((b_1, ..., b_i, ..., b_n))$$

$$b_k \in B_k \quad k = 1, ..., n \quad b_{i \cdot 1}, ..., b_{i \cdot p} \in B_i \quad \omega \in \Omega_2(p) \quad a \in A$$

- Докажем следующую лемму.

  *Лемма* 4.4. Для любого $c \in A$ эндоморфизм $f(c)$ $\Omega_2$-алгебры $M$ согласовано с эквивалентностью $N$.

- Пусть $\omega \in \Omega_2(p)$. Из равенства (4.3) следует

(4.4) $$f(c)((b_1, ..., b_{i \cdot 1}...b_{i \cdot p}\omega, ..., b_n)) = (b_1, ..., f_i(c)(b_{i \cdot 1}...b_{i \cdot p}\omega), ..., b_n)$$

Так как $f_i(c)$ - эндоморфизм $\Omega_2$-алгебры $B_i$, то из равенства (4.4) следует

(4.5) $$f(c)((b_1, ..., b_{i \cdot 1}...b_{i \cdot p}\omega, ..., b_n)) = (b_1, ..., f_i(c)(b_{i \cdot 1})...f_i(c)(b_{i \cdot p})\omega, ..., b_n)$$

Из равенств (4.5), (4.2) следует

(4.6) $$\begin{aligned}&f(c)((b_1, ..., b_{i \cdot 1}...b_{i \cdot p}\omega, ..., b_n))\\&= (b_1, ..., f_i(c)(b_{i \cdot 1}), ..., b_n)...(b_1, ..., f_i(c)(b_{i \cdot p}), ..., b_n)\omega\end{aligned}$$

Из равенств (4.6), (4.3) следует

(4.7) $$\begin{aligned}&f(c)((b_1, ..., b_{i \cdot 1}...b_{i \cdot p}\omega, ..., b_n))\\&= f(c)((b_1, ..., b_{i \cdot 1}, ..., b_n))...f(c)((b_1, ..., b_{i \cdot p}, ..., b_n))\omega\end{aligned}$$

Так как $f_i(c)$ - эндоморфизм $\Omega_2$-алгебры $B_i$, то из равенства (4.7) следует

(4.8) $$\begin{aligned}&f(c)((b_1, ..., b_{i \cdot 1}...b_{i \cdot p}\omega, ..., b_n))\\&= f(c)((b_1, ..., b_{i \cdot 1}, ..., b_n))...(b_1, ..., b_{i \cdot p}, ..., b_n)\omega\end{aligned}$$

Следовательно, мы доказали следующее утверждение.

*Лемма* 4.5. Равенство (4.8) следует из равенства (4.2).

- Следующее утверждение следует из равенства (4.3).

  *Лемма* 4.6. Равенство

(4.9) $$f(c)((b_1, ..., f_i(a)(b_i), ..., b_n)) = f(c)(f(a)((b_1, ..., b_i, ..., b_n)))$$

следует из равенства (4.3).

- Из лемм 4.5, 4.6 и определения [3]-2.2.13 следует лемма 4.4.

Из леммы 4.4 и теоремы [3]-2.2.14 следует, что на множестве $^*M/N$ определена
$\Omega_1$-алгебра. Рассмотрим диаграмму

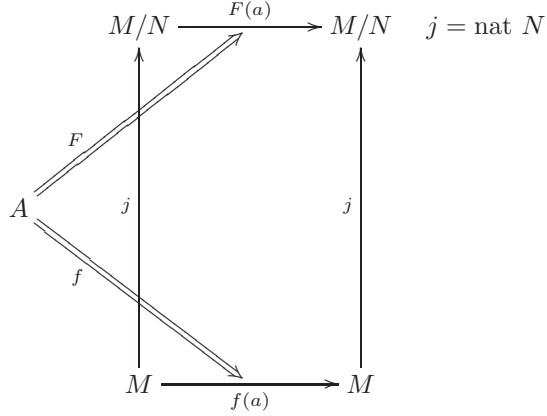

Согласно лемме 4.4, из условия

$$j(b_1) = j(b_2)$$

следует

$$j(f(a)(b_1)) = j(f(a)(b_2))$$

Следовательно, преобразование $F(a)$ определено корректно и

(4.10) $$F(a) \circ j = j \circ f(a)$$

Если $\omega \in \Omega_1(p)$, то мы положим

$$(F(a_1)...F(a_p)\omega)(J(b)) = J((f(a_1)...f(a_p)\omega)(b))$$

Следовательно, отображение $F$ является представлением $\Omega_1$-алгебры $A$. Из (4.10) следует, что $(\mathrm{id}, j)$ является морфизмом представлений $f$ и $F$.

Рассмотрим коммутативную диаграмму

(4.11)

$$\begin{array}{c} B_1 \times ... \times B_n \xrightarrow{i} M \xrightarrow{j} M/N \\ \searrow^{g_1} \nearrow \end{array}$$

Из равенств (4.1), (4.2), (4.3) следует

(4.12) $$g_1((b_1, ..., b_{i\cdot 1}...b_{i\cdot p}\omega, ..., b_n))$$
$$= g_1((b_1, ..., b_{i\cdot 1}, ..., b_n))...g_1((b_1, ..., b_{i\cdot p}, ..., b_n))\omega$$

(4.13) $$g_1((b_1, ..., f_i(a)(b_i), ..., b_n))$$
$$= f(a)(g_1((b_1, ..., b_i, ..., b_n)))$$

Из равенств (4.12) и (4.13) следует, что отображение $g_1$ является приведённым полиморфизмом представлений $f_1, ..., f_n$.

Поскольку $B_1 \times ... \times B_n$ - базис представления $M$ $\Omega_1$-алгебры $A$, то, согласно теореме [3]-3.2.7, для любого представления

$$A \xrightarrow{*} V$$

и любого приведенного полиморфизма

$$g_2 : B_1 \times ... \times B_n \longrightarrow V$$

существует единственный морфизм представлений $k : M \to V$, для которого коммутативна следующая диаграмма

(4.14)
$$\begin{array}{c} B_1 \times ... \times B_n \xrightarrow{i} M \\ \searrow_{g_2} \quad \downarrow_{k} \\ V \end{array}$$

Так как $g_2$ - приведенный полиморфизм, то $\ker k \supseteq N$.

Согласно теореме 3.3 отображение $j$ универсально в категории морфизмов представления $f$, ядро которых содержит $N$. Следовательно, определён морфизм представлений

$$h : M/N \to V$$

для которого коммутативна диаграмма

(4.15)
$$\begin{array}{c} M/N \\ \nearrow_{j} \quad \downarrow_{h} \\ M \searrow_{k} \\ V \end{array}$$

Объединяя диаграммы (4.11), (4.14), (4.15), получим коммутативную диаграмму

$$\begin{array}{c} B_1 \times ... \times B_n \xrightarrow{i} M \xrightarrow{j} M/N \\ \searrow_{g_2} \quad \downarrow_{k} \quad \downarrow_{h} \\ V \end{array}$$

Так как $\operatorname{Im} g_1$ порождает $M/N$, то отображение $h$ однозначно определено. $\square$

Согласно доказательству теоремы 4.3

$$B_1 \otimes ... \otimes B_n = M/N$$

Для $d_i \in A_i$ будем записывать

(4.16) $\qquad j \circ (d_1, ..., d_n) = d_1 \otimes ... \otimes d_n$

**Теорема 4.7.** *Пусть $B_1, ..., B_n$ - $\Omega_2$-алгебры. Пусть*

$$f : B_1 \times ... \times B_n \to B_1 \otimes ... \otimes B_n$$

*приведенный полиморфизм, определённый равенством*

(4.17) $\qquad f \circ (b_1, ..., b_n) = b_1 \otimes ... \otimes b_n$

*Пусть*

$$g : B_1 \times ... \times B_n \to V$$

приведенный полиморфизм в $\Omega_2$-алгебру $V$. Существует морфизм представлений
$$h : B_1 \otimes ... \otimes B_n \to V$$
такой, что диаграмма

$$\begin{array}{c}
B_1 \times ... \times B_n \xrightarrow{f} B_1 \otimes ... \otimes B_n \\
\searrow_{g} \qquad \downarrow_{h} \\
V
\end{array}$$

коммутативна.

*Доказательство.* Равенство (4.17) следует из равенств (4.1) и (4.16). Существование отображения $h$ следует из определения 4.1 и построений, выполненных при доказательстве теоремы 4.3. □

Равенства (4.12) и (4.13) можно записать в виде

$$(4.18) \quad \begin{aligned} &b_1 \otimes ... \otimes (b_{i \cdot 1}...b_{i \cdot p}\omega) \otimes ... \otimes b_n \\ &= (b_1 \otimes ... \otimes b_{i \cdot 1} \otimes ... \otimes b_n)...(b_1 \otimes ... \otimes b_{i \cdot p} \otimes ... \otimes b_n)\omega \end{aligned}$$

$$(4.19) \quad b_1 \otimes ... \otimes (f_i(a)(b_i)) \otimes ... \otimes b_n = f(a)(b_1 \otimes ... \otimes b_i \otimes ... \otimes b_n)$$

$$b_k \in B_k \quad k = 1, ..., n \quad b_{i \cdot 1}, ..., b_{i \cdot p} \in B_i \quad \omega \in \Omega_2(p) \quad a \in A$$

## 5. Список литературы

# 6. Предметный указатель

# 7. Специальные символы и обозначения

$B^{\otimes n}$  тензорная степень представления 8

$B_1 \otimes ... \otimes B_n$  тензорное произведение представлений 8